# ADAPTIVE ESTIMATION OF AND ORACLE INEQUALITIES FOR PROBABILITY DENSITIES AND CHARACTERISTIC FUNCTIONS[1]

By Sam Efromovich

*University of Texas at Dallas*

The theory of adaptive estimation and oracle inequalities for the case of Gaussian-shift–finite-interval experiments has made significant progress in recent years. In particular, sharp-minimax adaptive estimators and exact exponential-type oracle inequalities have been suggested for a vast set of functions including analytic and Sobolev with any positive index as well as for Efromovich–Pinsker and Stein blockwise-shrinkage estimators. Is it possible to obtain similar results for a more interesting applied problem of density estimation and/or the dual problem of characteristic function estimation? The answer is "yes." In particular, the obtained results include exact exponential-type oracle inequalities which allow to consider, for the first time in the literature, a simultaneous sharp-minimax estimation of Sobolev densities with any positive index (not necessarily larger than 1/2), infinitely differentiable densities (including analytic, entire and stable), as well as of not absolutely integrable characteristic functions. The same adaptive estimator is also rate minimax over a familiar class of distributions with bounded spectrum where the density and the characteristic function can be estimated with the parametric rate.

**1. Introduction.** Univariate probability density estimation is one of the fundamental topics in applied and mathematical statistics, and it is not surprising that first theoretical results about rate-optimal estimation of nonparametric functions were obtained for this statistical model; the interested reader is referred to a discussion in books [9, 14, 45, 49, 51]. An important step in the theory of a nonparametric density estimation was made by Nussbaum [42] who established that, for the case of a finite-support density and

Received February 2005; revised April 2007.
[1]Supported in part by NSF Grants DMS-02-43606, DMS-06-04558 and NSA Grant MDA904-07-1-0075.

a bounded loss function, there existed an asymptotic equivalence between the density model and a Gaussian-shift–finite-interval experiment; the interested reader can find more about the equivalence and a review of latest results in [5]. Because a Gaussian-shift model is simpler to work with, over the last decade the nonparametric research has been primarily devoted to a Gaussian-shift experiment and a vast set of pioneering results, specifically in the area of adaptive estimation and oracle inequalities, has been obtained; see a discussion in [6, 8, 15, 21, 24, 37, 39, 44, 54].

Due to Nussbaum's equivalence paradigm, there is a belief in the nonparametric literature that known adaptive estimators and oracle inequalities for a Gaussian-shift–finite-interval experiment may guide a creation of similar results for density estimation. This article shows that this belief is valid, and it develops a theory of adaptive estimation of and oracle inequalities for the probability density which matches recently obtained results for Gaussian-shift models. Moreover, it is possible to consider densities with both finite and infinite supports while the equivalence theory exists only for the density with a finite support, and the article also explores estimation of characteristic functions.

There are many applications of the obtained results. In particular, exponential-type oracle inequalities allow the statistician to consider a vast portfolio of blocks and thresholds including the smaller blocks suggested in the Gaussian-shift literature. The article also solves a long (more than two decades) standing problem of adaptive-sharp-minimax estimation of densities with a positive Sobolev index. Let us recall that, under mean integrated squared error (MISE) criteria, so far only densities with Sobolev index larger than 1/2 have been studied in the sharp-minimax literature; see a discussion in [3, 16, 18, 19, 23, 29, 32, 46, 47, 48, 50]. Note that, according to [17], no such restriction exists for a Gaussian-shift experiment. Interestingly, the asymptotic nonequivalence between the two models is valid whenever the index is at most 1/2, and for years this fact has served as a pleasing justification of the absence of the theory of a sharp adaptive estimation of those rougher densities; see a discussion in [4, 19]. This article shows that, fortunately, the nonequivalence does not affect the studied adaptive density estimation under the MISE criteria. Another important application is the possibility to consider distributions with not absolutely integrable (but square-integrable) characteristic functions which never before have been studied in the literature, and then suggest oracle inequalities for and sharp-minimax estimators of such characteristic functions. Further, for the first time in the literature a data-driven procedure for estimation of densities supported on a real line is suggested which is simultaneously sharp minimax over Sobolev (of any order) and infinitely differentiable densities (including entire densities like normal and their mixtures or analytic densities like Cauchy and their mixtures). Moreover, the suggested estimator implies the parametric rate of convergence for classical distributions with bounded spectrum (whose Fourier transform has a finite support).

The content of the article is as follows. To make the paper shorter, results are presented for densities supported on a real line (technical report [22] contains results for the finite support). Section 2 presents a short review of relevant results for the case of a Gaussian-shift experiment; these are the results to match. Section 3 presents the EP estimators for density and characteristic functions. Section 4 presents new oracle inequalities. Section 5 explores minimaxity of the estimator. The Stein density estimator, based on the famous Stein shrinkage procedure, is explored in Section 6; it is shown that, under a mild assumption, Stein and EP estimators have similar asymptotic properties. Discussion of results is deferred until Section 7. Section 8 contains proofs; some of its technically involved parts, including new moment and exponential inequalities for Sobolev statistics, are placed in the Appendix.

In what follows $C$'s denote generic positive constants and $o_s(1)$'s denote generic finite sequences which vanish as $s \to \infty$.

## 2. Review of relevant results for a Gaussian-shift experiment.

Consider a Gaussian-shift–finite-interval experiment $dY(t) = f(t) + n^{-1/2} dB(t)$, $0 \leq t \leq 1$, where $Y(t)$ is an observed signal, $f$ is an unknown square-integrable signal/shift of interest, $B(t)$ is a standard Brownian motion and $n$ is a positive integer which later will denote the sample size in a density model. Note that another customarily used name for the problem is the filtering a signal from a white Gaussian noise. Traditionally the model is rewritten in Fourier, wavelet or any other orthogonal basis domain; then an equivalent sequence model is considered:

$$y_j = \theta_j + n^{-1/2}\xi_j, \qquad j = 1, 2, \ldots, \tag{2.1}$$

where $\xi_k$ are independent standard Gaussian random variables, $\theta = \{\theta_1, \theta_2, \ldots\}$ is an unknown vector-parameter of interest, and $\int_0^1 f^2(t)\, dt = \sum_{j=1}^\infty \theta_j^2 =: \|\theta\|^2 < \infty$. The interested reader can find a comprehensive discussion of the sequence model (2.1) in [36]. The Efromovich–Pinsker (EP) blockwise-shrinkage estimator is defined as

$$\tilde{\theta}_j := \sum_{k=1}^K \tilde{\mu}_k y_j I(j \in B_k), \tag{2.2}$$

where the shrinkage (smoothing) coefficients/weights are

$$\tilde{\mu}_k := \frac{\|y\|_k^2 - L_k n^{-1}}{\|y\|_k^2} I(\|y\|_k^2 \geq (1 + t_k) L_k n^{-1}), \tag{2.3}$$

$I(\cdot)$ is the indicator, $\{1 = b_1 < b_2 < \cdots\}$ is a given sequence of positive integers and then $B_k := \{b_k, b_k + 1, \ldots, b_{k+1} - 1\}$ and $L_k := b_{k+1} - b_k$ are corresponding blocks and their lengths, $t_k > 0$ are thresholds (some authors refer to $1 + t_k$ as a penalty), $\|y\|_k^2 := \sum_{j \in B_k} y_j^2$ and this statistic is often referred to as a Sobolev statistic, and an integer $K = K(n)$ is a cutoff defined from the relation $\sum_{k=1}^K L_k < n^{1-1/\ln(n+1)} \le \sum_{k=1}^{K+1} L_k$ (see a comment on this choice in Section 7). The risk $E\|\mu_k y - \theta\|_k^2$ is minimized by a shrinkage coefficient (oracle)

$$(2.4) \qquad \mu_k := \frac{\|\theta\|_k^2}{\|\theta\|_k^2 + L_k n^{-1}},$$

which depends on a quantity $\|\theta\|_k^2 := \sum_{j \in B_k} \theta_j^2$ (so-called Sobolev functional) unavailable to the statistician. Then $\tilde{\theta}^* := (\tilde{\theta}_1^*, \tilde{\theta}_2^*, \ldots)$ with $\tilde{\theta}_j^* := \mu_k y_j$, $j \in B_k$ can serve as a (linear) blockwise-shrinkage oracle which, in its turn, is a blockwise version of the famous Wiener filter. The oracle has excellent minimax properties; in particular under a mild assumption on blocks and thresholds this oracle is simultaneously sharp minimax over Sobolev and analytic function classes; the interested reader can find a discussion in [17, 19, 36, 52, 53]. Then it is natural to use the mean squared error [or mean integrated squared error (MISE) for the dual filtering problem] of this oracle as a benchmark for the risk of any blockwise-shrinkage estimator. A simple calculation yields that the oracle's risk is

$$
\begin{aligned}
(2.5) \qquad E\|\tilde{\theta}^* - \theta\|^2 &= \sum_{k=1}^K \sum_{j \in B_k} E(\mu_k y_j - \theta_j)^2 + \sum_{k>K} \|\theta\|_k^2 \\
&= n^{-1} \sum_{k=1}^K L_k \mu_k + \sum_{k>K} \|\theta\|_k^2.
\end{aligned}
$$

Now we can formulate a known technical result which will imply oracle inequalities of interest. To do this, let us recall the Stirling formula for the Gamma function $\Gamma(L/2)$, $L = 1, 2, \ldots$ (see [1]),

$$(2.6) \qquad 1 < s_L^* \le \frac{\Gamma(L/2)}{(2\pi)^{1/2} e^{-L/2} (L/2)^{(L/2)-1/2}} \le s_L^{**} < \infty,$$
$$s_L^{**} \to 1 \text{ as } L \to \infty.$$

LEMMA 2.1 ([21]).   *Consider a particular block $B_k$ and assume that $0 < t_k \le 1$. Then there exists an absolute constant $C_0$ such that for any $q_k \in [1/4, \min(1, 1/4t_k))$ and any $\nu_k > 0$ the EP estimator satisfies*

$$(2.7) \qquad E\|\tilde{\theta} - \theta\|_k^2 \le E\|\tilde{\theta}^* - \theta\|_k^2 + n^{-1} L_k[\mu_k D_k^* + D_k^{**}],$$

where $E\|\tilde{\theta}^* - \theta\|_k^2 = n^{-1}\mu_k L_k$,

$$
\begin{aligned}
(2.8) \quad D_k^* := \nu_k + (1 + \nu_k^{-1})[&C_0^{1/2} L_k^{-1}(1 + (1 - q_k^{1/2})^{-2} t_k^{-1}) \\
&+ C_0 \mu_k (L_k t_k^2)^{-2}(1 + 2t_k)^3 \\
&+ \min(\mu_k(1 + t_k), 2t_k) I(\|\theta\|_k^2 < 2L_k t_k n^{-1})]
\end{aligned}
$$

and

$$
\begin{aligned}
(2.9) \quad D_k^{**} := (1 + \nu_k^{-1}) L_k^{-1} [&L_k^{1/2}/s_{L_k}^* + 8((L_k t_k)^{-1/4} + (L_k t_k^2)^{-1/2})] \\
&\times \exp\{-L_k[q_k t_k - \ln(1 + q_k t_k)]/2\} \\
&\times I(\|\theta\|_k^2 < (1 - q_k^{1/2})^2 L_k t_k n^{-1}).
\end{aligned}
$$

REMARK 2.1. The condition $t_k \to 0$, $k \to \infty$ is necessary for the estimate to sharply mimic the oracle's risk; this explains why only the case $t_k \le 1$ is considered in Lemma 2.1; see [24]. Further, it is easy to recognize that in (2.7) the term $D_k^{**}$ plays more important role than the $D_k^*$; indeed $D_k^{**}$ defines the remainder in the oracle inequality while $D_k^*$ defines the multiplicative factor. As a result, for mimicking the oracle's risk the term $D_k^{**}$ should vanish with an appropriate rate while the term $D_k^*$ may vanish with any rate as $k \to \infty$. Further, the following lower bound of [21]:

$$
\begin{aligned}
(2.10) \quad E\|\tilde{\theta} - \theta\|_k^2 \ge \frac{t_k}{s_{L_k}^{**}(1 + t_k)} &n^{-1} L_k^{1/2} \\
&\times \exp\{-L_k[t_k - \ln(1 + t_k)]/2\}, \qquad \|\theta\|_k = 0,
\end{aligned}
$$

allows one to appreciate the accuracy of the exponential factor in $D_k^{**}$ [compare exponential factors in (2.9) and (2.10)]. The exponential factor in $D_k^{**}$ is critical because it allows one to use smaller blocks; see a discussion in [6, 7, 8, 12, 19].

Now we can formulate several types of oracle inequalities suggested in the literature and based on Lemma 2.1. These are the results to match for density estimation.

THEOREM 2.1 ([21]). *Suppose that the assumption of Lemma 2.1 holds for all $k \in \{1, 2, \ldots, K\}$. Then:*

(a) *Risk of the EP estimate is bounded from above by the following oracle inequality:*

$$
(2.11) \quad E\|\tilde{\theta} - \theta\|^2 \le E\|\tilde{\theta}^* - \theta\|^2 + n^{-1} \sum_{k=1}^{K} L_k[\mu_k D_k^* + D_k^{**}].
$$

(b) *Denote* $\Delta_m := \max_{m \le k \le K} D_k^*$, $S_m := \sum_{k=m}^{K} L_k D_k^{**}$ *and* $\Upsilon_0 := \{k : \mu_k D_k^* + D_k^{**} \ge 1\}$. *Then*

$$
\begin{aligned}
(2.12) \quad E\|\tilde{\theta} - \theta\|^2 \le \min_{1 \le m \le K} & \left[ (1 + \Delta_m) E\|\tilde{\theta}^* - \theta\|^2 + n^{-1}\left( S_m + \sum_{k=1}^{m-1} L_k \right) \right] \\
& + n^{-1} \sum_{k \in \Upsilon_0} L_k [\mu_k D_k^* + D_k^{**}],
\end{aligned}
$$

*where by convention* $\sum_{j=1}^{0} = 0$.

(c) *Set* $\bar{\Upsilon}_0 := \{k : \bar{D}_k^* + \bar{D}_k^{**} \ge 1\}$ *with* $\bar{D}_k^*$ *and* $\bar{D}_k^{**}$ *defined as in* (2.8) *and* (2.9) *only with* $\mu_k$ *and indicator functions replaced by 1, and then, following part* (b), *define corresponding* $\bar{\Delta}_m$, $\bar{S}_m$ *and* $\bar{\Upsilon}_0$. *Also, let us modify the EP estimator* $\hat{\theta}_j$ *by considering* $\check{\theta}_j := y_j$ *for* $j \in B_k$, $k \in \bar{\Upsilon}_0$ *and* $\check{\theta}_j := \hat{\theta}_j$ *otherwise. Then*

$$
\begin{aligned}
(2.13) \quad E\|\check{\theta} - \theta\|^2 \le \min_{1 \le m \le K} & \left[ (1 + \bar{\Delta}_m) E\|\tilde{\theta}^* - \theta\|^2 + n^{-1}\left( \bar{S}_m + \sum_{k=1}^{m-1} L_k \right) \right] \\
& + n^{-1} \sum_{k \in \bar{\Upsilon}_0} L_k.
\end{aligned}
$$

## 3. EP density and characteristic function estimators.

Suppose that $X_1, \ldots, X_n$, $n > 3$ are i.i.d. realizations according to an unknown square-integrable on a real line density $f(x)$, $x \in (-\infty, \infty)$; it is not assumed that the density is positive on a real line. Let us recall that

$$
(3.1) \qquad f(x) = (2\pi)^{-1} \int_{-\infty}^{\infty} h(u) e^{-iux} \, du, \qquad x \in (-\infty, \infty),
$$

where

$$
(3.2) \qquad h(u) := E\{e^{iuX}\} = \int_{-\infty}^{\infty} e^{iux} f(x) \, dx, \qquad u \in (-\infty, \infty)
$$

is the characteristic function corresponding to $f$. If the characteristic function is not absolutely integrable, then the inverse formula (3.1) is understood in the sense of Plancherel's theorem. The problem is to estimate the density and the characteristic function under the MISE criterion.

Recall that the characteristic function satisfies $h(-u) = \overline{h(u)}$, the complex conjugate of $h(u)$. Thus we can consider only $h(u)$, $u \in [0, \infty)$ and then $f(x) = \pi^{-1} \int_0^{\infty} \mathrm{Re}\{h(u) e^{-iux}\} \, du$. Now we are following the construction of the EP estimator for the Gaussian-shift case. We divide a half-line $[0, \infty)$ into a sequence of nonoverlapping blocks (intervals) $B_k := [b_k', b_{k+1}')$, $0 = b_1' < b_2' < \cdots$ with the corresponding lengths $L_k := b_{k+1}' - b_k' = \int_{B_k} du$. Then the following abuse of the previous notation will be handy. Set

$$
(3.3) \qquad \|y\|_k^2 := \int_{B_k} |\hat{h}(u)|^2 \, du,
$$

where

$$(3.4) \qquad \hat{h}(u) := n^{-1} \sum_{l=1}^{n} \exp\{iuX_l\}$$

is the empirical characteristic function estimator. Then we define an EP density estimator as

$$(3.5) \qquad \tilde{f}(x) := \pi^{-1} \int_0^\infty \mathrm{Re}\{\tilde{h}(u)e^{-iux}\}\, du, \qquad x \in (-\infty, \infty).$$

Here

$$(3.6) \qquad \tilde{h}(u) := \sum_{k=1}^{K} \tilde{\mu}_k \hat{h}(u) I(u \in B_k), \qquad u \geq 0$$

is the EP characteristic function estimator, and $\tilde{\mu}_k$ is defined in (2.3). To make the similarity complete, we denote

$$(3.7) \qquad \|\theta\|_k^2 := \int_{B_k} |h(u)|^2\, du, \qquad \|\tilde{\theta}\|_k^2 := \int_{B_k} |\tilde{h}(u)|^2\, du$$

and

$$(3.8) \qquad \|\tilde{\theta} - \theta\|_k^2 := \int_{B_k} |\tilde{h}(u) - h(u)|^2\, du.$$

To shed light on the above-introduced notation, note that according to Plancherel's identity the MISE of EP density estimator (3.5) can be written as

$$(3.9) \qquad \begin{aligned} E \int_{-\infty}^{\infty} (\tilde{f}(x) - f(x))^2\, dx &= \pi^{-1} E \int_0^\infty |\tilde{h}(u) - h(u)|^2\, du \\ &= \pi^{-1} E \sum_{k=1}^{\infty} \|\tilde{\theta} - \theta\|_k^2. \end{aligned}$$

Further, using (2.4) we define the corresponding oracles $\tilde{f}^*(x)$ and $\tilde{h}^*(u)$ as

$$(3.10) \qquad \tilde{f}^*(x) := \pi^{-1} \sum_{k=1}^{K} \mu_k \int_{B_k} \mathrm{Re}\{\hat{h}(u)e^{-iux}\}\, du, \qquad x \in (-\infty, \infty),$$

$$(3.11) \qquad \tilde{h}^*(u) := \sum_{k=1}^{K} \mu_k \hat{h}(u) I(u \in B_k), \qquad u \geq 0.$$

Also we set $\|\tilde{\theta}^*\|_k^2 := \int_{B_k} |\tilde{h}^*(u)|^2\, du$ and $\|\tilde{\theta}^* - \theta\|_k^2 := \int_{B_k} |\tilde{h}^*(u) - h(u)|^2\, du$.

Finally, if an EP density estimate (or the oracle) takes on negative values, then its nonnegative projection may be considered; see Section 3.1 in [19]. Further, if a monotonicity assumption is known, then methods of Efromovich [20] can be used.

REMARK 3.1. According to (3.9), for both the density and characteristic function settings, it suffices to present bounds on their MISEs via $E\|\tilde\theta - \theta\|^2 := E\sum_{k=1}^\infty \|\tilde\theta - \theta\|_k^2$. This approach will be used in Section 4. Also, to avoid any confusion with the Gaussian-shift case, we shall refer to the above-introduced $\tilde\theta$ as the EP density-model estimator.

**4. Exponential-type oracle inequality for EP estimator.** In what follows $c_1$ denotes the universal positive constant $C_2$ of de la Peña and Montgomery-Smith [13], $c_2$ denotes the universal positive constant $K$ in the Bernstein-type inequality (3.18) of Giné, Latala and Zinn [27], $d := d(f) := \int_{-\infty}^{\infty} |h(u)|^2\,du = 2\pi\int_{-\infty}^{\infty} f^2(x)\,dx$, $d^* := d^*(f, L) := \min_{z>0}(z + Lz^{-1}\int_{\{x:\,f(x)\geq z\}} f^2(x)\,dx)$, and for a $k$th block

$$
\begin{aligned}
\lambda_1 &:= \lambda_1(L_k, t_k, d, d^*) \\
&:= (dc_1^2 c_2)^{-1}(1 - \min(1/2, t_k^{1/4}))^2(1 - (L_k + 1)^{-1/2})^2 \\
&\quad \times \min\Bigg(\frac{1}{[1 + 4n^{-1}t_k(2d^{-1/2} + 3n^{-1}d^{-1}t_k)]}, \\
&\qquad\qquad \frac{c_1 d}{t_k[8d^*(f, L_k) + 3(n^{-1}L_k t_k)^{1/2}]}, \\
&\qquad\qquad \frac{d[nc_1^4 t_k^{-4} L_k^{-5/2}]^{1/3}}{[(2d)^{1/2} + 20n^{-1}t_k]^{1/3}}, \frac{c_1^{3/2} dn^{1/2}}{2t_k^{3/2} L_k}\Bigg),
\end{aligned}
\tag{4.1}
$$

$$
\begin{aligned}
\lambda_2 &:= \lambda_2(L_k, t_k, d) \\
&:= \frac{n\min(1/4, t_k^{1/2})}{L_k t_k c_1^2}\frac{(1 - \min(1/2, t_k^{1/4}))^2(1 - (L_k + 1)^{-1/2})^2}{3c_1^{-1} + 2dL_k^{-1}t_k^{-1} + 8n^{-1}(2d^{1/2} + t_k)},
\end{aligned}
\tag{4.2}
$$

$$
\lambda_3 := \lambda_3(L_k, t_k, d) := \frac{\min(1/4, t_k^{1/2})}{6t_k d^{1/2}}\frac{(1 - (L_k + 1)^{-1/2})^2}{1 + (n^{-1}t_k L_k^{3/2}d^{-1})^{1/2}}.
\tag{4.3}
$$

THEOREM 4.1. *Suppose that $X_1, \ldots, X_n$, $n > 3$ are i.i.d. according to a square-integrable density $f \in L_2(-\infty, \infty)$. Consider a particular block $B_k$ with length $L_k > 0$ and a particular threshold level $t_k > 0$. Then for any $\nu_k \in (0, 1)$ the following oracle inequality holds for the EP density-model estimator defined in Section 3:*

$$
E\|\tilde\theta - \theta\|_k^2 \leq E\|\tilde\theta^* - \theta\|_k^2 + n^{-1}L_k[\mu_k D_k' + D_k''],
\tag{4.4}
$$

*where*

$$
E\|\tilde\theta^* - \theta\|_k^2 = n^{-1}L_k\mu_k[1 - \mu_k L_k^{-1}\|\theta\|_k^2],
\tag{4.5}
$$

$$
D_k' := \nu_k(1 - \mu_k L_k^{-1}\|\theta\|_k^2) + (1 + \nu_k^{-1})
$$

$$(4.6) \qquad \times [L_k^{-1/2}(15d^{1/2} + 3d(1 + L_k^{-1/2})(1 + t_k^{-1}))$$
$$+ \min(\mu_k(1 + t_k), 2t_k)I(\|\theta\|_k^2 < 2L_k t_k n^{-1})],$$

$$(4.7) \qquad D_k'' := (1 + \nu_k^{-1})[L_k^{-1}(d + 3d^{1/2}t_k)]^{1/2}$$
$$\times G(L_k, t_k, d, d^*)I(\|\theta\|_k^2 < L_k^{1/2}t_k n^{-1}),$$

$$(4.8) \qquad G(L_k, t_k, d, d^*) := [c_1 c_2 \exp\{-t_k^2 L_k \lambda_1\}$$
$$+ 2c_1 \exp\{-t_k^2 L_k \lambda_2\} + \exp\{-t_k^2 L_k \lambda_3\}]^{1/2}.$$

Theorem 4.1 implies a result which matches Theorem 2.1.

COROLLARY 4.1. *Let the assumption of Theorem 4.1 hold for all $k \in \{1, 2, \ldots, K\}$. Then assertions* (a)–(c) *of Theorem 2.1 are valid for the EP density-model estimator with $D_k^*$ and $D_k^{**}$ replaced by $D_k'$ and $D_k''$, respectively.*

These results yield two important conclusions: (i) It is possible to suggest identical blockwise-shrinkage estimators for the Gaussian-shift, density and characteristic function estimation models. (ii) The MISEs of those data-driven estimators satisfy similar exponential-type oracle inequalities.

REMARK 4.1. While there is a difference between the density-model oracle's error $E\|\tilde{\theta}^* - \theta\|_k^2$, presented in (4.5), and the corresponding Gaussian-shift oracle's error $E\|\tilde{\theta}^* - \theta\|_k^2 = n^{-1}L_k\mu_k$, this difference bears no consequences for nonparametric cases where the MISE vanishes more slowly than $n^{-1}$. The latter is based on a plain observation that for any square-integrable density the term $\mu_k L_k^{-1}\|\theta\|_k^2$ vanishes as $k \to \infty$; further, note that $\mu_k L_k^{-1}\|\theta\|_k^2 \leq L_k^{-1}\pi d$ and if the statistician uses blocks satisfying $L_k \geq L(n) \to \infty$, $n \to \infty$, then this term vanishes uniformly over the blocks as $n \to \infty$.

REMARK 4.2. In a majority of asymptotic applications of Theorem 4.1 the main exponential term in (4.8) is the one containing $\lambda_1$. Further, let us note that $d^*(f, L) \leq 2\min(\sup_x f(x), (dL)^{1/2})$. This inequality allows one to analyze $\lambda_1$ for bounded and unbounded densities.

## 5. Sharp minimaxity.

In this section the above-established oracle inequality is used to prove a simultaneous sharp minimaxity of the EP density estimate for Sobolev and infinitely differentiable distribution classes as well as its rate minimaxity for distribution classes with bounded spectrum where the MISE converges with the parametric rate $n^{-1}$. The interested reader can find a thorough discussion of these classes in [3, 19, 30, 31, 32, 33, 34, 35, 36, 38, 40, 51, 55]. Below these distribution classes are defined via corresponding characteristic functions which are assumed to be square integrable, and let us recall that if the characteristic function $h$ belongs to $L_2(-\infty, \infty)$, then the corresponding cumulative distribution function is absolutely continuous and its density $f$ belongs to $L_a(-\infty, \infty)$ for any $1 \le a \le 2$; see Theorem 11.6.1 in [38].

We consider those distribution classes in turn. Let $\alpha$ and $Q$ be positive real numbers; then a Sobolev class (of order $\alpha$) is defined as

$$
(5.1) \qquad \mathcal{S}(\alpha, Q) := \Big\{ f(x) : \pi^{-1} \int_0^\infty (1 + |u|^{2\alpha}) |h(u)|^2 \, du \le Q,
$$
$$
h(u) = \int_{-\infty}^\infty f(x) e^{iux} \, dx \Big\}.
$$

THEOREM 5.1 (Sobolev class). *Let a sample $X_1, X_2, \ldots, X_n$ of $n$ i.i.d. observations with a square-integrable density $f \in L_2(-\infty, \infty)$ be given. Suppose that blocks and thresholds of EP estimator $\tilde{f}$, defined in (3.5), satisfy*

$$
(5.2) \qquad L_{k+1}/L_k \to 1 \quad and \quad \sup_f D_k' \to 0
$$
$$
as \ k \to \infty, \ \sup_f \sum_{k=1}^K L_k D_k'' < \delta_n,
$$

*where the supremums are taken over $f \in \mathcal{S}(\alpha, Q)$ and $\delta_n = n^{o_n(1)}$. Then*

$$
(5.3) \qquad \sup_{f \in \mathcal{S}(\alpha, Q)} \Big\{ E \int_{-\infty}^\infty (\tilde{f}(x) - f(x))^2 \, dx \Big\} (1 + o_n(1))
$$

$$
(5.4) \qquad = \inf_{\check{f}} \sup_{f \in \mathcal{S}(\alpha, Q)} E \int_\infty^\infty (\check{f}(x) - f(x))^2 \, dx
$$
$$
= P(\alpha, Q) n^{-2\alpha/(2\alpha+1)} (1 + o_n(1)),
$$

*where in (5.4) the infimum is taken over all possible density estimates $\check{f}$ based on the sample and parameters $\alpha$ and $Q$, and $P(\alpha, Q) := (2\alpha+1)[\pi(2\alpha+1)(\alpha+1)\alpha^{-1}]^{-2\alpha/(2\alpha+1)} Q^{1/(2\alpha+1)}$ is the Pinsker constant.*

Let us recall that only Sobolev classes of order $\alpha > 1/2$ have been considered in the literature so far; see a discussion in [11, 19, 29, 46, 50].

Now let us consider another popular (specifically in the literature devoted to characteristic functions and stable distributions) class of infinitely differentiable distributions

$$
(5.5) \qquad \mathcal{A}(r, \gamma, Q)
$$
$$
:= \Big\{ f : \pi^{-1} \int_0^\infty |e^{\gamma u^r} h(u)|^2 \, du \le Q, \ h(u) = \int_{-\infty}^\infty f(x) e^{iux} \, dx \Big\}.
$$

Here $\gamma$ and $Q$ are positive real numbers and $r \in (0, 2]$. A thorough discussion of this class can be found in the classical books [38, 40, 55] as well as in [2, 19, 33, 34, 35, 39]. This class includes analytic, stable and entire distributions with more familiar particular examples being Cauchy mixtures (where $r = 1$) and Normal mixtures (where $r = 2$).

THEOREM 5.2 (Infinitely differentiable class). *Let the assumption of Theorem 5.1 hold with the supremums in (5.2) taken over $f \in \mathcal{A}(r, \gamma, Q)$ and $\delta_n = o_n(1)[\ln(n)]^{1/2}$. Then*

$$(5.6) \qquad \sup_{f \in \mathcal{A}(r, \gamma, Q)} \left\{ E \int_{-\infty}^{\infty} (\tilde{f}(x) - f(x))^2 \, dx \right\} (1 + o_n(1))$$

$$= \inf_{\check{f}} \sup_{f \in \mathcal{A}(r, \gamma, Q)} E \int_{-\infty}^{\infty} (\check{f}(x) - f(x))^2 \, dx$$

$$(5.7) \qquad = \pi^{-1} n^{-1} [\ln(n)/(2\gamma)]^{1/r} (1 + o_n(1)),$$

*where the infimum in (5.7) is taken over all estimates $\check{f}$ based on the sample and parameters $(r, \gamma, Q)$.*

Finally, let $s$ denote a positive real number, and consider a familiar class of distributions with bounded spectrum

$$(5.8) \qquad \mathcal{B}(s) = \left\{ f : h(u) = 0, \ |u| > s, \ h(u) = \int_{-\infty}^{\infty} f(x) e^{iux} \, dx \right\}.$$

According to Theorem 11.12.1 in [38], a distribution with bounded spectrum, which is not from a uniform family, is an entire order of 1 and of exponential type. Then, as it is emphasized by Ibragimov and Khasminskii [33, 34], we are dealing with essentially infinite-dimensional class. Nonetheless, they were the first to recognize that the sharp-minimax MISE is $\pi^{-1} s n^{-1}(1 + o_n(1))$, that is, the MISE's convergence is *parametric*! The parametric convergence is too fast for the essentially nonparametric adaptive EP estimator; however, the following result still holds.

THEOREM 5.3 (Bounded spectrum class). *Let the assumption of Theorem 5.1 hold with the supremums in (5.2) taken over $f \in \mathcal{B}(s)$ and $\delta_n = o_s(1)s$. Then*

$$(5.9) \qquad \sup_{f \in \mathcal{B}(s)} \left\{ E \int_{-\infty}^{\infty} (\tilde{f}(x) - f(x))^2 \, dx \right\} (1 + o_n(1) + o_s(1))$$

$$(5.10) \qquad = \inf_{\check{f}} \sup_{f \in \mathcal{B}(s)} E \int_{-\infty}^{\infty} (\check{f}(x) - f(x))^2 \, dx = \pi^{-1} s n^{-1}(1 + o_n(1)),$$

*where the infimum in (5.10) is taken over all estimates $\check{f}$ based on the sample and parameter $s$.*

REMARK 5.1. Using Remark 4.2, it is plain to verify that a majority of known portfolios of blocks and thresholds, suggested in the sharp-minimax Gaussian-shift literature, simultaneously satisfy conditions of Theorems 5.1–5.3. Just to point to a specific and simple example with relatively "small" logarithmic blocks, consider $\{(L_k = \ln^3(k+3),\ t_k = 1/\ln(\ln(k+3))), k = 1, 2, \ldots\}$. This portfolio simultaneously satisfies assumptions of Theorems 5.1–5.3.

We may conclude that the adaptive EP density (or characteristic function) estimator is *simultaneously* sharp minimax over Sobolev and infinitely differentiable classes of distributions. On top of this nice property, the adaptive estimator is also rate minimax over distributions with bounded spectrum, and its MISE attains the parametric-minimax MISE when the spectrum band increases. To the best of the author's knowledge, this is the first known example of such simultaneous adaptive density estimation, as well as the first example of a simultaneous adaptive sharp-minimax estimation for classes of distributions which include both absolutely integrable and not absolutely integrable characteristic functions.

REMARK 5.2. Let us note that due to Plancherel's identity, results of Theorems 5.1–5.3, except for using an extra factor $2\pi$ in the formulas for minimax MISEs, hold for the dual problem of characteristic function estimation. As a result, the EP characteristic function estimator (3.6) is simultaneously sharp minimax over Sobolev and infinitely differentiable distribution classes, and it is also rate minimax over classes of distributions with bounded spectrum.

**6. Stein estimator.** The blockwise-shrinkage literature, devoted to Gaussian-shift experiments, also explores a Stein (blockwise-shrinkage) estimator which, using notation of Section 2, can be written as

$$(6.1) \quad \bar{\theta}_j := \frac{\|y\|_k^2 - (1 + t_k) L_k n^{-1}}{\|y\|_k^2} I(\|y\|_k^2 \geq (1 + t_k) L_k n^{-1}) y_j, \qquad j \in B_k.$$

Note that if the EP estimator uses a hard block-thresholding, a Stein estimator uses a soft one. Then, according to the paradigm of Section 3, the Stein density estimator can be defined as

$$(6.2) \qquad \bar{f}_S(x) := \pi^{-1} \int_0^\infty \mathrm{Re}\{\bar{h}_S(u) e^{-iux}\}\, du,$$

where the Stein characteristic function estimator is

$$(6.3) \qquad \bar{h}_S(u) := \sum_{k=1}^K \bar{\mu}_k \hat{h}(u) I(u \in B_k), \qquad u \geq 0$$

and, recalling notation (3.3),

$$(6.4) \qquad \bar{\mu}_k := \frac{\|y\|_k^2 - (1+t_k)L_k n^{-1}}{\|y\|_k^2} I(\|y\|_k^2 \geq (1+t_k)L_k n^{-1}).$$

The following proposition allows one to explore the Stein density (or characteristic function) estimator via its EP counterpart. Recall that $G(L, t, d, d^*)$ was defined in (4.8).

THEOREM 6.1. *Let* $\tilde{f}$ *and* $\bar{f}_S$ *denote EP and Stein estimators which use the same blocks, thresholds and* $K$. *Suppose that the assumption of Theorem 4.1 holds. Then*

$$E \int_{-\infty}^{\infty} (\bar{f}_S(x) - \tilde{f}(x))^2 \, dx$$

$$\leq \pi^{-1} n^{-1} \sum_{k=1}^{K} L_k \mu_k$$

$$(6.5) \qquad \times [12 L_k^{-1/2}(1 - (L_k+1)^{-1/2})^{-2}(d^{1/2} + dt_k^{-1}(1 + L_k^{-1/2}))$$

$$+ 2t_k I(\|\theta\|_k^2 \geq (1/2)L_k t_k n^{-1})]$$

$$+ \pi^{-1} n^{-1} \sum_{k=1}^{K} L_k t_k^2 (1+t_k)^{-1} G^2(L_k, t_k/2, d, d^*)$$

$$\times I(\|\theta\|_k^2 < (1/2)L_k^{1/2} t_k n^{-1}).$$

This result implies that, under the MISE criteria and for the portfolios of blocks and thresholds discussed in Section 5, the two estimators perform similarly.

## 7. Discussion.

7.1. *Parameter* $K$ *in EP estimator.* In the theory of oracle inequalities this parameter is assumed to be given; see [8, 21]. For a minimax (or adaptive) setting it should be chosen in such a way that the squared bias of the oracle is negligible with respect to its variance. For instance, for a Sobolev class with index $\alpha$ this is achieved if $\tilde{h}^*(u)$ is zero on frequencies larger than $\gamma_n n^{-1/(2\alpha+1)}$ where $\gamma_n$ increases to infinity as slowly as desired; this remark explains how $K := K(n)$ was chosen in Section 2. At the same time, for infinitely differentiable distributions $K(n)$ may be logarithmic, that is, dramatically smaller than for Sobolev functions. Further, for distributions with bounded spectrum $K(n)$ may be any increasing-to-infinity sequence.

7.2. *Sobolev classes with index $\alpha \leq 1/2$.* This is a new addition to the set of distributions covered by the theory of minimax estimation and oracle inequalities. Obviously there were serious technical difficulties in dealing with such Sobolev densities. Also, the case of Sobolev densities with index larger than $1/2$ is very nice and appealing because it implies that the characteristic function is absolutely integrable, the corresponding density is defined by Fourier inverse formula (3.1) and it is bounded and uniformly continuous. Sobolev characteristic functions with index $\alpha \leq 1/2$ do not have these nice properties and, moreover, the inverse formula (3.1) is understood, according to Plancherel's theorem, as a limit in $L_2(-\infty, \infty)$-norm of $(2\pi)^{-1} \int_{-A}^{A} e^{-iux} h(u) \, du$, $A \to \infty$. At the same time, it is important to note that the characteristic function is not necessarily absolutely integrable and, for instance, there is a vast class of Pólya-type characteristic functions that are not absolutely integrable. Namely, according to the famous Pólya condition, a real-valued and continuous function $g(u)$, $u \in (-\infty, \infty)$ is the characteristic function of an absolutely continuous distribution if $g(0) = 1$, $g(-u) = g(u)$, $g(u)$ is convex for positive $u$ and $g(u) \to 0$ as $u \to \infty$; see [40], page 70. Note that the condition involves no restriction on how fast $g(u)$ must vanish. Characteristic functions $h(u) = [1 + |u|^{\beta}]^{-1}$, $1/2 < \beta \leq 1$ from the Linnik distribution family as well as $h(u) = [1 + |u|^2]^{-\rho}$, $1/4 < \rho \leq 1/2$, studied by Karl Pearson, are particular examples of characteristic functions which are square integrable but not absolutely integrable; see [2, 43].

7.3. *Why MISE?* A choice of the loss function in the density estimation literature has been always a source of hot debates thanks to statisticians passionately devoted to $L_1$-distance, different $L_p$-distances with $p > 1$, Hellinger distances, distances based upon Kullback–Leibler numbers, etc. The interested reader can find a discussion of these approaches in [9, 14, 19]. Until now, there was no objective argument in favor of the $L_2$-distance/MISE because it was always assumed that underlying characteristic functions were absolutely integrable. The inclusion of not absolutely integrable characteristic functions changes the situation because now Plancherel's theorem (and correspondingly $L_2$-norm) is the necessary tool. This remark, at least partially, may serve as a justification for using the MISE criteria.

7.4. *Distributions with finite support.* In many applied problems the statistician knows support of the density; circular data is a familiar example. Suppose that the support is $[0, 1]$. Then, following the Gaussian-shift approach of Section 2, the density can be written in Fourier domain as $f(x) = [1 + \sum_{j=1}^{\infty} \theta_j \varphi_j(x)] I(x \in [0, 1])$, $\theta_j := \int_0^1 f(x) \varphi_j(x)$ where $\{1, \varphi_j(x) := 2^{1/2} \cos(\pi j x), j = 1, 2, \ldots\}$ is a classical cosine basis on $[0, 1]$. Further, $\|y\|_k^2 := \sum_{j \in B_k} y_j^2$, $y_j := n^{-1} \sum_{l=1}^{n} \varphi_j(X_l)$ may serve as an analogue of $\|y\|_k^2$ in the Gaussian-shift and probability settings of Sections 2 and 3 (note that here we again intentionally use the same notation). Then the EP finite-support density estimator is defined as

$$(7.1) \qquad \tilde{f}(x) := \left[ 1 + \sum_{k=1}^{K} \tilde{\mu}_k \sum_{j \in B_k} y_j \varphi_j(x) \right] I(x \in [0,1]),$$

where $K$ and $\tilde{\mu}_k$ are the same as in Sections 2 and 3. Corresponding exponential inequalities and minimax results can be found in the technical report [22].

7.5. *Different types of oracle inequalities.* Corollary 4.1 (or Theorem 2.1) presents three different types of oracle inequalities, and each may be useful on its own. Inequality (2.13) is useful because all its components, apart from the oracle's MISE, depend only on blocks and thresholds but not on an estimated function. This type of oracle inequalities can be found in [6, 7, 8]. The other types of inequalities, originated in [18], are more complicated because the remainder depends on an estimated function; but this complexity may be useful. As an example, let us present a discussion of the, phenomenon, mentioned in Remark 2.1, of the necessity for thresholds to vanish for sharp mimicking of the oracle's MISE. The nonparametric blockwise-shrinkage literature contains results of intensive numerical studies which indicate excellent performance of estimates with nonvanishing thresholds; see a discussion in [6, 7, 10, 12]. Do these studies contradict the theory? To answer this question, let us examine oracle inequalities in Theorem 2.1. Oracle inequality (2.13) cannot shed light on the phenomenon because $t_k$ must vanish for the right-hand side of (2.13) to converge to the oracle's MISE. On the other hand, oracle inequalities (2.11) and (2.12) can explain the phenomenon. Indeed, we can relax the assumption on thresholds by assuming that an estimated density (or shift function) satisfies

$$\frac{\sum_{k=1}^{K} L_k \mu_k \min(\mu_k(1+t_k), 2t_k) I(\|\theta\|_k^2 < 2 L_k t_k n^{-1})}{\sum_{k=1}^{K} L_k \mu_k} = o_n(1).$$

It is not difficult to check numerically that this assumption often holds for functions used in numerical studies. Thus, oracle inequalities (2.11) and (2.12) have allowed us to shed a new light on the above-mentioned numerical results.

Further, let us note that an assumption like (5.2), by including terms depending on an estimated density, bears the same flavor as the oracle inequalities (2.11) and (2.12) because it allows the statistician to justify/explain a special portfolio of blocks and thresholds for a targeted class of functions.

7.6. *Possible applications in related problems.* There are many related applied problems where the obtained results may motivate new research and innovative procedures, with particular examples being survival analysis, deconvolution, biased data, error density estimation, time series analysis, etc. The developed methodology can be of special interest for the analysis of wavelet estimators. A discussion of possible extensions can be found in [5, 22, 25, 26, 28, 41, 45].

**8. Proofs.** In what follows $\hat{\Theta}_k := L_k^{-1}\|y\|_k^2 - n^{-1}$ and $\Theta_k := L_k^{-1}\|\theta\|_k^2$.

PROOF OF THEOREM 4.1. A direct calculation implies that

$$(8.1) \qquad E\hat{h}(u) = h(u), \qquad E|\hat{h}(u) - h(u)|^2 = n^{-1}(1 - |h(u)|^2).$$

Recall that $E\|\tilde{\theta}^* - \theta\|_k^2 = E\int_{B_k}|\mu_k\hat{h}(u) - h(u)|^2\,du$ and write for any $\nu \in (0,1)$,

$$
\begin{aligned}
(8.2) \qquad E\|\tilde{\theta} - \theta\|_k^2 &= E\int_{B_k}|\tilde{h}(u) - h(u)|^2\,du \\
&= E\int_{B_k}|\tilde{\mu}_k\hat{h}(u) - h(u)|^2\,du \\
&= E\int_{B_k}|(\mu_k\hat{h}(u) - h(u)) + (\bar{\mu}_k - \mu_k)\hat{h}(u)|^2\,du \\
&\leq (1+\nu)E\int_{B_k}|\mu_k\hat{h}(u) - h(u)|^2\,du \\
&\quad + (1+\nu^{-1})E\int_{B_k}|(\tilde{\mu}_k - \mu_k)\hat{h}(u)|^2\,du.
\end{aligned}
$$

Using (8.1) we get

$$
\begin{aligned}
(8.3) \qquad & E\int_{B_k}|\mu_k\hat{h}(u) - h(u)|^2\,du \\
&= E\int_{B_k}|\mu_k(\hat{h}(u) - h(u)) - (1 - \mu_k)h(u)|^2\,du \\
&= \mu_k^2 E\int_{B_k}|\hat{h}(u) - h(u)|^2\,du + (1 - \mu_k)^2\int_{B_k}|h(u)|^2\,du \\
&= n^{-1}\mu_k^2\int_{B_k}(1 - |h(u)|^2)\,du + (1 - \mu_k)^2 L_k\Theta_k \\
&= L_k\left[\frac{\Theta_k^2 n^{-1}}{(\Theta_k + n^{-1})^2} + \frac{n^{-2}\Theta_k}{(\Theta_k + n^{-1})^2}\right] - n^{-1}\mu_k^2 L_k\Theta_k \\
&= n^{-1}L_k\mu_k - n^{-1}\mu_k^2 L_k\Theta_k.
\end{aligned}
$$

In particular, this verifies (4.5). The second expectation in the right-hand side of (8.2) can be written as

$$E\left\{\int_{B_k} |(\tilde{\mu}_k - \mu_k)\hat{h}(u)|^2 \, du\right\} = E\{(\tilde{\mu}_k - \mu_k)^2 L_k(\hat{\Theta}_k + n^{-1})\} =: E\{A\}.$$

Let us evaluate the term $A$. Note that (2.3) can be rewritten as $\tilde{\mu}_k = \hat{\Theta}_k(\hat{\Theta}_k + n^{-1})^{-1} I(\hat{\Theta}_k \geq t_k n^{-1})$. In what follows we skip subscripts whenever no confusion may occur. Write

$$A = \frac{n^{-2}(\hat{\Theta} - \Theta)^2 I(\hat{\Theta} \geq tn^{-1})L}{(\hat{\Theta} + n^{-1})(\Theta + n^{-1})^2} + \mu^2 L(\hat{\Theta} + n^{-1})I(\hat{\Theta} < tn^{-1}) =: A_1 + A_2.$$

Set $q := 1 - (L+1)^{-1/2}$ and evaluate $A_1$:

$$A_1 = \frac{n^{-2}(\hat{\Theta} - \Theta)^2 I(\hat{\Theta} \geq tn^{-1})L}{(\hat{\Theta} + n^{-1})(\Theta + n^{-1})^2}[I(\Theta < (1-q)tn^{-1}) + I(\Theta \geq (1-q)tn^{-1})]$$

$$\leq \frac{n^{-2}(\hat{\Theta} - \Theta)^2 L}{(\hat{\Theta} + n^{-1})(\Theta + n^{-1})^2}I(\hat{\Theta} - \Theta > qtn^{-1})I(\Theta < (1-q)tn^{-1})$$

$$+ \frac{n^{-1}(\hat{\Theta} - \Theta)^2 L}{(\Theta + n^{-1})^2}I(\Theta \geq (1-q)tn^{-1}) =: A_{11} + A_{12}.$$

Plainly $A_{11} \leq L(\hat{\Theta} - \Theta)I(\hat{\Theta} - \Theta > qtn^{-1})I(\Theta < (1-q)tn^{-1})$, and using the Cauchy–Schwarz inequality we get

$$E\{A_{11}\} \leq LE^{1/2}\{(\hat{\Theta} - \Theta)^2\}\mathrm{Pr}^{1/2}\{\hat{\Theta} - \Theta > qtn^{-1}\}I(\Theta < (1-q)tn^{-1}).$$

To continue we need a result that will be proved in the Appendix.

LEMMA 8.1. *Let the assumption of Theorem* 4.1 *hold. Set* $d := \int_{-\infty}^{\infty} |h(u)|^2 \, du$, $d_j := \max_{v \in B} \int_B (|h(u-v)|^j + |h(u+v)|^j) \, du$, $d^* := \min_{z>0}(z + Lz^{-1}\int_{\{x: f(x) \geq z\}} f^2(x) \, dx)$. *Then:*

(a) *The moment inequality holds:*

$$(8.4) \qquad E(\hat{\Theta} - \Theta)^2 \leq L^{-1}n^{-1}[2d_1\Theta + d_2 n^{-1}].$$

(b) *For* $q = 1 - (L+1)^{-1/2}$,

$$(8.5) \qquad \mathrm{Pr}\{\hat{\Theta} - \Theta > qtn^{-1}\}I(\Theta < (1-q)tn^{-1}) \leq G^2(L, t, d, d^*),$$

*where* $G(L, t, d, d^*)$ *is defined in* (4.8).

(c) *The following relations between* $d_1$, $d_2$ *and* $d$ *hold:*

$$(8.6) \qquad d_1 \leq [2Ld_2]^{1/2} \quad and \quad d_2 \leq d.$$

Using Lemma 8.1 we get

$$
(8.7) \quad
\begin{aligned}
&E\{A_{11}\} \\
&\leq n^{-1} L [L^{-1}(d + 2^{3/2} d^{1/2} t)]^{1/2} G(L, t, d, d^*) I(\Theta < (1-q) t n^{-1}).
\end{aligned}
$$

Note that $(1-q)^{-1} = (L+1)^{1/2}$, and then using (8.4) and (8.6) we get

$$
\begin{aligned}
E\{A_{12}\} &\leq n^{-1} L (\Theta + n^{-1})^{-2} L^{-1} n^{-1} [2d_1 \Theta + d_2 n^{-1}] I(\Theta \geq (1-q) t n^{-1}) \\
&\leq n^{-1} L [L^{-1}(\mu 2 d_1 + d_2)] I(\Theta \geq (1-q) t n^{-1}) \\
&\leq n^{-1} \mu L [L^{-1}(2^{3/2} L^{1/2} d^{1/2} + d[1 + (L+1)^{1/2} t^{-1}])] \\
&\leq n^{-1} \mu L [L^{-1/2}(2^{3/2} d^{1/2} + d(L^{-1/2} + (1 + L^{-1/2}) t^{-1}))].
\end{aligned}
$$

Further,

$$
\begin{aligned}
A_2 &= \mu^2 L (\hat{\Theta} + n^{-1}) I(\hat{\Theta} < t n^{-1}) [I(\Theta \geq 2 t n^{-1}) + I(\Theta < 2 t n^{-1})] \\
&\leq \mu^2 L n^{-1} (1+t) I(\Theta - \hat{\Theta} \geq \Theta/2) I(\Theta \geq 2 t n^{-1}) \\
&\quad + \mu^2 L n^{-1} (1+t) I(\Theta < 2 t n^{-1}) =: A_{21} + A_{22}.
\end{aligned}
$$

Using the Chebyshev inequality and (8.4) we get

$$
\begin{aligned}
E\{A_{21}\} &\leq n^{-1} \mu L \left[ (1+t) \mu \frac{n^{-1} L^{-1}(2 d_1 \Theta + d_2 n^{-1})}{(\Theta/2)^2} I(\Theta \geq 2 t n^{-1}) \right] \\
&\leq n^{-1} \mu L [2 L^{-1}(4(2 L d_2)^{1/2} + d_2 t^{-1})] \\
&\leq n^{-1} \mu L [12 d^{1/2} L^{-1/2} + 2 d L^{-1} t^{-1}].
\end{aligned}
$$

To evaluate $A_{22}$ we note that $(1+t)\mu I(\Theta < 2 t n^{-1}) \leq 2 t I(\Theta < 2 t n^{-1})$, and then

$$
A_{22} \leq n^{-1} \mu L [\min(\mu(1+t), 2t) I(\Theta < 2 t n^{-1})].
$$

Combining the obtained results we conclude that

$$
(8.8) \quad
\begin{aligned}
E\|\tilde{\theta} - \theta\|_k^2 &\leq E\|\tilde{\theta}^* - \theta\|_k^2 \\
&\quad + n^{-1} L_k \mu_k [\nu_k (1 - \mu_k \Theta_k + (1 + \nu_k^{-1})) \\
&\qquad \times [L_k^{-1/2}(15 d^{1/2} + 3d(1 + L_k^{-1/2})(1 + t_k^{-1})) \\
&\qquad\quad + \min(\mu_k(1 + t_k), 2 t_k) I(\Theta_k < 2 t_k n^{-1})]] \\
&\quad + n^{-1} L_k (1 + \nu_k^{-1}) [L_k^{-1}(d + 3 d^{1/2} t_k)]^{1/2} \\
&\qquad \times G(L_k, t_k, d, d^*) I(\Theta_k < L_k^{-1/2} t_k n^{-1}).
\end{aligned}
$$

This verifies (4.4).  $\square$

PROOF OF THEOREM 5.1. Relation (5.4) is known for the case of $\alpha > 1/2$; see [23, 46, 50]. Consider the case $\alpha \leq 1/2$. The lower minimax bound (5.4) is established in [22]; the proof is too lengthy to reproduce it here and the interested reader is referred to [22]. The upper minimax bound (5.4), as well as the validity of upper bound (5.3) for any $\alpha$, is established with the help of the oracle inequality. We begin with the analysis of oracle $\tilde{f}^*$ defined in (3.10). Following [17], a direct calculation, based on using (4.5), yields that whenever $L_{k+1}/L_k \to 1$, $k \to \infty$ the oracle's MISE satisfies

$$\sup_{f \in \mathcal{S}(\alpha, Q)} E \int_{-\infty}^{\infty} (\tilde{f}^*(x) - f(x))^2 \, dx$$
$$= P(\alpha, Q) n^{-2\alpha/(2\alpha+1)} (1 + o_n(1)).$$

In other words, the oracle is sharp minimax.

Then oracle inequality (4.4), together with assumption (5.2), yields

$$\sup_{f \in \mathcal{S}(\alpha, Q)} E \int_{-\infty}^{\infty} (\tilde{f}(x) - f(x))^2 \, dx (1 + o_n(1))$$
$$= \sup_{f \in \mathcal{S}(\alpha, Q)} E \int_{-\infty}^{\infty} (\tilde{f}^*(x) - f(x))^2 \, dx$$
$$= P(\alpha, Q) n^{-2\alpha/(2\alpha+1)} (1 + o_n(1)).$$

This result shows that for Sobolev classes the EP estimator is sharp minimax and matches performance of the oracle. $\square$

Theorems 5.2 and 5.3 are verified identically.

PROOF OF THEOREM 6.1. Write

$$\int_{-\infty}^{\infty} (\bar{f}_S(x) - \tilde{f}(x))^2 \, dx$$
$$= \pi^{-1} \int_0^{\infty} |\bar{h}_S(u) - \tilde{h}(u)|^2 \, du = \pi^{-1} \sum_{k=1}^{K} (\bar{\mu}_k - \tilde{\mu}_k)^2 \|y\|_k^2$$
(8.9)
$$= \pi^{-1} \sum_{k=1}^{K} \frac{(L_k t_k n^{-1})^2}{\|y\|_k^2} I(\|y\|_k^2 \geq (1 + t_k) L_k n^{-1})$$
$$\leq \pi^{-1} n^{-1} \sum_{k=1}^{K} L_k t_k^2 (1 + t_k)^{-1} I(\hat{\Theta}_k \geq t_k n^{-1})$$
$$\times [I(\Theta_k < (t_k/2) n^{-1}) + I(\Theta_k \geq (t_k/2) n^{-1})].$$

Now let us make several preliminary calculations. First of all, set $q :=$ $1 - (1 + L_k)^{-1/2}$, and write

$$I(\hat{\Theta}_k \geq t_k n^{-1}) I(\Theta_k < (t_k/2) n^{-1})$$
$$\leq I(\hat{\Theta}_k - \Theta_k > (t_k/2) n^{-1}) I(\Theta_k < (t_k/2) n^{-1})$$
$$\leq I(\hat{\Theta}_k - \Theta_k > q(t_k/2) n^{-1})$$
$$\times [I(\Theta_k < (1-q)(t_k/2) n^{-1})$$
$$+ I((1-q)(t_k/2) n^{-1} \leq \Theta_k < (t_k/2) n^{-1})].$$

Using (8.5) we get

$$E\{I(\hat{\Theta}_k - \Theta_k > q(t_k/2) n^{-1})\} I(\Theta_k < (1-q)(t_k/2) n^{-1})$$
$$\leq G^2(L_k, t_k/2, d, d^*) I(\Theta_k < (1 + L_k)^{-1/2}(t_k/2) n^{-1}).$$

Using (8.4), (8.6) and a plain inequality $I(\Theta_k \geq bn^{-1}) \leq \mu_k (1 + b^{-1})$ we get

$$E\{I(\hat{\Theta}_k - \Theta_k > q(t_k/2) n^{-1})\}$$
$$\times I((1-q)(t_k/2) n^{-1} \leq \Theta < (t_k/2) n^{-1})$$
$$\leq 4 L_k^{-1} n^{-1} [2 d_1 \Theta_k + d_2 n^{-1}] [q^2 t_k^2 n^{-2}]^{-1}$$
$$\times I((1-q)(t_k/2) n^{-1} \leq \Theta_k < (t_k/2) n^{-1})$$
$$\leq 4 q^{-2} t_k^{-2} L_k^{-1} [2 \mu_k (\Theta_k + n^{-1}) n (2 L_k d)^{1/2}$$
$$+ d \mu_k (1 + 2(1 + L_k)^{1/2} t_k^{-1})] I(\Theta_k < (t_k/2) n^{-1})$$
$$\leq 4 q^{-2} t_k^{-2} \mu_k L_k^{-1/2} [3(1 + t_k) d^{1/2} + d(L_k^{-1/2} + 2 t_k^{-1}(1 + L_k^{-1/2}))]$$
$$\times I(\Theta_k < (t_k/2) n^{-1}).$$

Combining these results we get

$$\pi^{-1} n^{-1} \sum_{k=1}^{K} L_k t_k^2 (1 + t_k)^{-1} E\{I(\hat{\Theta}_k \geq t_k n^{-1})\} I(\Theta_k < (t_k/2) n^{-1})$$
$$\leq \pi^{-1} n^{-1} \sum_{k=1}^{K} L_k t_k^2 (1 + t_k)^{-1} G^2(L_k, t_k/2, d, d^*)$$
$$\times I(\Theta_k < (L_k + 1)^{-1/2}(t_k/2) n^{-1})$$
$$+ \pi^{-1} n^{-1} \sum_{k=1}^{K} L_k \mu_k [12 L_k^{-1/2} (1 - (L_k + 1)^{-1/2})^{-2}$$
$$\times (d^{1/2} + d t_k^{-1}(1 + L_k^{-1/2}))].$$

Further, we note that $I(\Theta_k \geq (t_k/2)n^{-1}) \leq \mu_k(1+2t_k^{-1})$, and that

$$\pi^{-1}n^{-1}\sum_{k=1}^{K}L_k t_k^2(1+t_k)^{-1}I(\Theta_k \geq (t_k/2)n^{-1})$$

$$\leq \pi^{-1}n^{-1}\sum_{k=1}^{K}L_k\mu_k[2t_k I(\Theta_k \geq (t_k/2)n^{-1})].$$

Combining results verifies (6.5).  □

## APPENDIX

PROOF OF LEMMA 8.1.  The first part of (8.6) is based on the Cauchy–Schwarz inequality, the second on the remark that $-(u+v) \leq u-v$ for any $u,v \in B$ and $|h(u)| = |h(-u)|$ (let us also note that $d_2 \to d$ as $L \to \infty$). Let us now verify (8.4). Write

$$
\begin{aligned}
E(\hat{\Theta}-\Theta)^2 &= L^{-2}E\left[\int_B(|\hat{h}(u)|^2-|h(u)|^2-n^{-1})\,du\right]^2 \\
&= L^{-2}E\int_B\int_B(|\hat{h}(u)|^2-|h(u)|^2-n^{-1}) \\
&\qquad\qquad \times\,(|\hat{h}(v)|^2-|h(v)|^2-n^{-1})\,du\,dv \\
&= L^{-2}\int_B\int_B E\{|\hat{h}(u)|^2|\hat{h}(v)|^2\}\,du\,dv \\
&\quad -2L^{-2}\int_B E\{|\hat{h}(u)|^2\}\,du\int_B(|h(v)|^2+n^{-1})\,dv \\
&\quad +L^{-2}\left[\int_B(|h(u)|^2+n^{-1})\,du\right]^2 =: A_1+A_2+A_3.
\end{aligned}
$$

(A.1)

Consider these three addends in turn. In what follows $l_k \neq l_m \neq \cdots \neq l_q$ means that all these parameters are different, and recall the assumption $n > 3$. Write

$$
\begin{aligned}
&n^4 E\{|\hat{h}(u)|^2|\hat{h}(v)|^2\} \\
&\quad = \sum_{l_1,l_2,l_3,l_4=1}^{n} E\{\exp(iuX_{l_1}-iuX_{l_2}+ivX_{l_3}-ivX_{l_4})\} \\
&\quad \leq \sum_{l_1\neq l_2\neq l_3\neq l_4=1}^{n}|h(u)|^2|h(v)|^2 \\
&\qquad +\sum_{l_1\neq l_2\neq l_3=l_4=1}^{n}|h(u)|^2+\sum_{l_1=l_2\neq l_3\neq l_4=1}^{n}|h(v)|^2
\end{aligned}
$$

$$+ 2 \sum_{l_1 \neq l_2 = l_3 \neq l_4 = 1}^{n} |h(u)||h(v)|(|h(u-v)| + |h(u+v)|)$$

$$+ 2 \sum_{l_1 \neq l_2 = l_3 = l_4 = 1}^{n} |h(u)|^2$$

$$+ 2 \sum_{l_1 = l_2 = l_3 \neq l_4 = 1}^{n} |h(v)|^2 + \sum_{l_1 = l_4 \neq l_2 = l_3 = 1}^{n} (|h(u-v)|^2 + |h(u+v)|^2)$$

$$+ \sum_{l_1 = l_2 \neq l_3 = l_4 = 1}^{n} 1 + \sum_{l_1 = l_2 = l_3 = l_4 = 1}^{n} 1.$$

Then, using $2|h(u)||h(v)| \leq |h(u)|^2 + |h(v)|^2$, $|h(u)| \leq 1$, $(n-1)(n-2)(n-3) = n^3 - 6n^2 + 11n - 6$, $(n-1)(n-2) = n^2 - 3n + 2$ and simple algebra, we get

$$A_1 \leq \Theta^2 [1 - 6n^{-1} + 11n^{-2} - 6n^{-3}] + 2n^{-1}\Theta[1 - 3n^{-1} + 2n^{-2}]$$
$$+ 2L^{-1}n^{-1}d_1\Theta + 4n^{-2}\Theta[1 - n^{-1}] + n^{-2}L^{-1}d_2 + n^{-2}$$
$$= \Theta^2 - 2n^{-1}\Theta^2 + \Theta^2[-4n^{-1} + 11n^{-12} - 6n^{-3}]$$
$$+ 2n^{-1}\Theta - 2n^{-2}\Theta + L^{-1}n^{-1}[2d_1\Theta + d_2n^{-1}] + n^{-2}.$$

Further, (8.1) implies $E|\hat{h}(u)|^2 = |h(u)|^2 + n^{-1}(1 - |h(u)|^2)$, and we get

$$A_2 = -2L^{-2} \int_B (|h(u)|^2 + n^{-1}(1 - |h(u)|^2)) \, du \int_B (|h(v)|^2 + n^{-1}) \, dv$$
$$= -2[\Theta + n^{-1}(1 - \Theta)][\Theta + n^{-1}]$$
$$= -2\Theta^2 - 4n^{-1}\Theta - 2n^{-2} + 2n^{-1}\Theta^2 + 2n^{-2}\Theta.$$

Also $A_3 = (n^{-1} + \Theta)^2 = n^{-2} + 2n^{-1}\Theta + \Theta^2$. Combining the results in (A.1) and using $-4n^{-1} + 11n^{-2} - 6n^{-3} \leq 0$ for $n > 1$ we verify (8.4).

Let us check (8.5). Write

$$L\hat{\Theta} = \int_B |\hat{h}(u)|^2 \, du - Ln^{-1} = n^{-2} \sum_{1 \leq l,m \leq n} \int_B \exp\{iu(X_l - X_m)\} \, du - Ln^{-1}$$
$$= \left[ n^{-2} \sum_{l=1}^{n} \int_B du - Ln^{-1} \right] + n^{-2} \sum_{1 \leq l \neq m \leq n} \int_B \exp\{iu(X_l - X_m)\} \, du$$
$$= n^{-2}2 \sum_{1 \leq l < m \leq n} \int_B \cos(u(X_l - X_m)) \, du =: n^{-2}2 \sum_{1 \leq l < m \leq n} g(X_l - X_m).$$

Note that $g(x,y) := g(x - y)$ is a symmetric function in $(x,y)$ which can be viewed as a kernel of $U$-statistics. Thus we can use known exponential inequalities for $U$-statistics to analyze $\hat{\Theta}$. In what follows $X, X_1, \ldots, X_n, Y, Y_1, \ldots,$ $Y_n$ are i.i.d. random variables according to an underlying density $f$. Using Hoeffding's decomposition we continue:

$$
\begin{aligned}
L\hat{\Theta} = 2n^{-2} &\sum_{1 \le l < m \le n} H(X_l, X_m) \\
\text{(A.2)} \qquad &+ 2(n-1)n^{-2} \sum_{l=1}^{n} (E\{g(X_l - Y)|X_l\} - E\{g(X-Y)\}) \\
&+ (n-1)n^{-1} E\{g(X-Y)\} =: \tilde{A}_1 + \tilde{A}_2 + \tilde{A}_3,
\end{aligned}
$$

where $H(X,Y) := g(X-Y) - E\{g(X-Y)|X\} - E\{g(X-Y)|Y\} + E\{g(X-Y)\}$. A direct calculation shows that

$$
\text{(A.3)} \quad E\{g(X-Y)\} = \text{Re}\left\{ \int_B E e^{iu(X-Y)} \, du \right\} = \int_B |h(u)|^2 \, du = L\Theta
$$

and

$$
\text{(A.4)} \qquad E\{g(X-Y)|Y\} = \text{Re}\left\{ \int_B e^{iuY} h(-u) \, du \right\}.
$$

This implies

$$
\begin{aligned}
\text{(A.5)} \qquad H(X,Y) = &\int_B \cos(u(X-Y)) \, du \\
&- \text{Re}\left\{ \int_B (e^{iuX} + e^{iuY}) h(-u) \, du \right\} + L\Theta.
\end{aligned}
$$

According to Theorem 1 in [13], for all $z > 0$ there exists a universal constant $c_1$ such that

$$
\text{(A.6)} \qquad \Pr\{|\tilde{A}_1| > z\} \le c_1 \Pr\{|\tilde{A}_1^*| > z/c_1\},
$$

where $\tilde{A}_1^* := n^{-2} \sum_{1 \le l \ne m \le n} H(X_l, Y_m) = n^{-2} \sum_{1 \le l, m \le n} H(X_l, Y_m) - n^{-2} \sum_{l=1}^{n} H(X_l, Y_l)$ is a decoupled version of $\tilde{A}_1$. Using (A.2)–(A.6) we write for any $q, \gamma \in (0, 1)$,

$$
\begin{aligned}
\text{(A.7)} \qquad &\Pr\{\hat{\Theta} - \Theta > qtn^{-1}\} \\
&= \Pr\{\tilde{A}_1 + \tilde{A}_2 + [(n-1)/n]L\Theta - L\Theta > qtLn^{-1}\} \\
&\le \Pr\{\tilde{A}_1 + \tilde{A}_2 > qtLn^{-1}\} \le \Pr\{\tilde{A}_1 > \gamma qtLn^{-1}\} \\
&\quad + \Pr\{\tilde{A}_2 > (1-\gamma)qtLn^{-1}\} \\
&\le c_1 \Pr\{|\tilde{A}_1^*| > \gamma qtLn^{-1}/c_1\} + \Pr\{\tilde{A}_2 > (1-\gamma)qtLn^{-1}\} \\
&\le c_1 \Pr\left\{ n^{-2} \left| \sum_{1 \le l, m \le n} H(X_l, Y_m) \right| > \gamma^2 qtLn^{-1}/c_1 \right\}
\end{aligned}
$$

$$+ c_1 \Pr\left\{ n^{-2} \left| \sum_{l=1}^{n} H(X_l, Y_l) \right| > \gamma(1-\gamma)qtLn^{-1}/c_1 \right\}$$

$$+ \Pr\{\tilde{A}_2 > (1-\gamma)qtLn^{-1}\}.$$

Consider the first probability. $H(x,y)$ is symmetric in $(x,y)$ and it is a completely degenerated kernel in the sense that $E\{H(X,Y)|X\} = 0$. Thus we can use the following exponential inequality (3.18) of [27].

LEMMA A.1. *Let* $X, X_1, \ldots, X_n, Y, Y_1, \ldots, Y_n$ *be i.i.d. Consider a symmetric and completely degenerated kernel* $H(x,y)$. *Then there exists a universal constant* $c_2$ *such that for any* $z > 0$

$$\Pr\left\{ \left| \sum_{1 \le l, m \le n} H(X_l, Y_m) \right| > z \right\}$$

$$(A.8) \quad \le c_2 \exp\left\{ -\frac{1}{c_2} \min\left( \frac{z^2}{n^2 E\{H^2(X,Y)\}}, \frac{z}{n\|H\|_*}, \right. \right.$$

$$\left. \left. \frac{z^{2/3}}{[n\|E\{H^2(X,Y)|X\}\|_\infty]^{1/3}}, \frac{z^{1/2}}{\|H\|_\infty^{1/2}} \right) \right\},$$

*where*

$$\|H\|_* := \sup_{\psi_1, \psi_2} \{E\{H(X,Y)\psi_1(X)\psi_2(Y)\} : E\{\psi_1^2(X)\} \le 1, E\{\psi_2^2(Y)\} \le 1\},$$

$$\|E\{H^2(X,Y)|X\}\|_\infty := \sup_x E\{H^2(x,Y)\} \text{ and } \|H\|_\infty := \sup_{x,y} H(x,y).$$

Let us evaluate in turn the four components of the minimum in (A.8). Using (A.5) and $(a+b+c)^2 \le 2a^2 + 4(b^2+c^2)$ we get

$$E\{H^2(X,Y)\}$$

$$\le E\left\{ \int_B \int_B (\cos[(u-v)(X-Y)] + \cos[(u+v)(X-Y)])\, du\, dv \right\}$$

$$+ 4E\left| \int_B (e^{iuX} + e^{iuY})h(-u)\, du \right|^2 + 4(L\Theta)^2$$

$$= \int_B \int_B (|h(u-v)|^2 + |h(u+v)|^2)\, du\, dv$$

$$+ 4E \int_B \int_B (e^{iuX} + e^{iuY})(e^{-ivX} + e^{-ivY})h(-u)h(v)\, du\, dv + 4(L\Theta)^2$$

$$\le Ld_2 + 8 \int_B \int_B [h(u-v) + h(u)h(-v)]h(-u)h(v)\, du\, dv + 4(L\Theta)^2$$

$$\le Ld_2 + 8[d_2 L^3]^{1/2}\Theta + 12(L\Theta)^2.$$

In the last inequality we used $\int_B |h(u-v)| \, du \leq [Ld_2]^{1/2}$. Recall that we are considering only $\Theta < (1-q)tn^{-1}$, $q \in (0,1)$, and get

$$(A.9) \quad E\{H^2(X,Y)\} \leq L[d_2 + 4L^{1/2}(1-q)tn^{-1}(2d_2^{1/2} + 3L^{1/2}(1-q)tn^{-1})].$$

Now we are considering $\|H\|_*$. In what follows the supremum is taken over $\psi_1$ and $\psi_2$ such that $E\psi_j^2(X) \leq 1$, $j = 1, 2$. Write

$$\|H\|_* = \sup E\{H(X,Y)\psi_1(X)\psi_2(Y)\}$$

$$\leq \sup E\{g(X-Y)\psi_1(X)\psi_2(Y)\} + 2\sup E\{E\{g(X-Y)|X\}\psi_1(X)\} + L\Theta$$

$$=: D_1 + 2D_2 + L\Theta.$$

Introduce $A := \{x : f(x) < z\}$, $\gamma \in (0,1)$, and write

$$D_1 = \sup E\left\{\int_B (1/2)[e^{iu(X-Y)} + e^{-iu(X-Y)}]\psi_1(X)\psi_2(Y) \, du\right\}$$

$$= \sup \int_B |E\{e^{iuX}\psi_1(X)\}|^2 \, du$$

$$\leq (1+\gamma)\sup \int_B |E\{I(X \in A)e^{iuX}\psi_1(X)\}|^2 \, du$$

$$+ (1+\gamma^{-1})\sup \int_B |E\{I(X \in A^c)\psi_1(X)\}|^2 \, du$$

$$=: D_{11} + D_{12}.$$

Using the Plancherel identity we get

$$D_{11} \leq (1+\gamma)(2\pi)\sup \int_A f^2(x)\psi_1^2(x) \, dx \leq (1+\gamma)2\pi z.$$

Further, using the Cauchy–Schwarz inequality we get

$$D_{12} \leq (1+\gamma^{-1})L \int_{A^c} f(x) \, dx \sup \int_{-\infty}^{\infty} f(x)\psi_1^2(x) \, dx$$

$$\leq (1+\gamma^{-1})Lz^{-1} \int_{A^c} f^2(x) \, dx.$$

Set $\gamma = 0.2$ and get $D_1 \leq 8d^*$.

Further, $D_2 = \sup E\{\text{Re}\{\int_B e^{iuX}h(-u) \, du\}\psi_1(X)\} \leq \int_B |h(u)| \, du \leq 2L\Theta^{1/2}$. Plainly $\Theta \leq \min((1-q)tn^{-1}, 1)$, and this yields that

$$\|H\|_* \leq 8d^* + 2L\Theta^{1/2} + L\Theta \leq 8d^* + 3L(1-q)^{1/2}t^{1/2}n^{-1/2}.$$

Further, let us consider $\|E\{H^2(X,Y)|X\}\|_\infty$. Using $(a+b+c+d)^2 \leq 2a^2 + 4b^2 + 8(c^2 + d^2)$ we get

$$E\{H^2(x,Y)\} \leq 2E\{g^2(x-Y)\} + 4E^2\{g(x-Y)\}$$

$$+ 8E\{E^2\{g(X-Y)|Y\}\} + 8L^2\Theta^2$$

$$= 2E\left\{\left[\int_B \cos(u(x-Y))\,du\right]^2\right\} + 4\left[\mathrm{Re}\left\{\int_B e^{iux}h(-u)\,du\right\}\right]^2$$

$$+ 8E\left[\mathrm{Re}\left\{\int_B e^{iuY}h(-u)\,du\right\}\right]^2 + 8L^2\Theta^2$$

$$\leq (1/2)E\int_B\int_B (e^{iu(x-Y)} + e^{-iu(x-Y)})$$

$$\times (e^{iv(x-Y)} + e^{-iv(x-Y)})\,du\,dv + 20L^2\Theta$$

$$\leq \int_B\int_B (|h(u-v)| + |h(u+v)|)\,du\,dv + 20L^2\Theta$$

$$\leq L[d_1 + 20(1-q)tLn^{-1}].$$

Finally, we have a plain inequality $\sup_{x,y}|H(x,y)| \leq 4L$. Using these results and Lemma A.1 we get

$$(A.10) \quad \mathrm{Pr}\left\{\left|\sum_{1\leq l,m\leq n} H(X_l,Y_m)\right| > \gamma^2 qtnL/c_1\right\} \leq c_2\exp\left\{-\frac{\gamma^4 q^2 t^2 L}{dc_1^2 c_2}\nu_1\right\},$$

where

$$\nu_1 := \min\left(\frac{1}{[1 + 4L^{1/2}(1-q)tn^{-1}(2d^{-1/2} + 3d^{-1}L^{1/2}(1-q)tn^{-1})]},\right.$$

$$(A.11) \qquad \frac{c_1 d}{t[8d^* + 3L((1-q)t/n)^{1/2}]}, \frac{d[c_1^4 nt^{-4}L^{-2}]^{1/3}}{[d_1 + 20(1-q)Ltn^{-1}]^{1/3}},$$

$$\left.\frac{c_1^{3/2}dn^{1/2}}{2t^{3/2}L}\right).$$

To evaluate the second probability in (A.7), let us recall Bernstein's inequality.

LEMMA A.2. *Let* $Z_1,\dots,Z_n$ *be i.i.d.,* $|Z_1| < M$ *a.e.,* $E\{Z_1\} = 0$ *and* $\mathrm{Var}(Z_1) = \sigma^2 < \infty$. *Then for any* $z > 0$

$$(A.12) \qquad \begin{aligned} &\max\left(\mathrm{Pr}\left\{\sum_{l=1}^n Z_l < -z\right\}, \mathrm{Pr}\left\{\sum_{l=1}^n Z_l > z\right\}\right) \\ &\leq \exp\left\{-\frac{z^2}{2n\sigma^2 + (2/3)Mz}\right\}. \end{aligned}$$

This implies

$$\mathrm{Pr}\left\{\left|\sum_{l=1}^n H(X_l,Y_l)\right| > \gamma(1-\gamma)qtnL/c_1\right\}$$

$$\leq 2\exp\left\{-\frac{\gamma^2(1-\gamma)^2q^2t^2L}{n^{-1}Ltc_1^2[3c_1^{-1}+2dL^{-1}t^{-1}+8n^{-1}(2d^{1/2}+t)]}\right\}.$$

Let us consider the third probability in (A.7). Write

$$\tilde{A}_2 = \frac{2(n-1)}{n^2}\sum_{l=1}^{n}\left(\operatorname{Re}\left\{\int_B e^{iuX_l}h(-u)\,du\right\}-L\Theta\right) =: \frac{2(n-1)}{n^2}\sum_{l=1}^{n}V_l.$$

Plainly $|V_l| \leq \int_B |h(u)|\,du + L\Theta \leq 2L\Theta^{1/2}$. Also, $E\{V_l\}=0$ and

$$\operatorname{Var}(V_l) \leq E\int_B\int_B(1/4)[e^{iuX}h(-u)+e^{-iuX}h(u)]$$
$$\times [e^{ivX}h(-v)+e^{-ivX}h(v)]\,du\,dv$$
$$\leq (1/2)d_1L\Theta.$$

Then Lemma A.2, $n>3$ and $\Theta<(1-q)tn^{-1}$ imply that

$$\Pr\{\tilde{A}_2 > (1-\gamma)qtLn^{-1}\}$$
$$= \Pr\left\{\sum_{l=1}^{n}V_l > 2^{-1}(1-\gamma)qtLn(n-1)^{-1}\right\}$$
$$\leq \exp\left\{-\frac{2^{-2}(1-\gamma)^2q^2t^2L^2[n/(n-1)]^2}{nd_1L\Theta+(2/3)2L\Theta^{1/2}2^{-1}(1-\gamma)qtL[n/(n-1)]}\right\}$$
$$\leq \exp\left\{-\frac{(1-\gamma)^2q^2t^2L}{4[(1-q)d_1t+(1-q)^{1/2}t^{3/2}Ln^{-1/2}]}\right\}.$$

Combining the obtained inequalities in (A.7) implies [$\nu_1$ is defined in (A.11)]

$$(A.13)\quad\begin{aligned}&\Pr\{\hat{\Theta}-\Theta>qtn^{-1}\}I(\Theta<(1-q)tn^{-1})\\&\quad\leq c_1c_2\exp\{-\gamma^2q^2t^2L\nu_1/(dc_1^2c_2)\}\\&\quad+2c_1\exp\left\{-\frac{\gamma^2(1-\gamma)^2q^2t^2L}{n^{-1}Ltc_1^2[3c_1^{-1}+2dL^{-1}t^{-1}+8n^{-1}(2d^{1/2}+t)]}\right\}\\&\quad+\exp\left\{-\frac{(1-\gamma)^2q^2t^2L}{4[(1-q)d_1t+(1-q)^{1/2}t^{3/2}Ln^{-1/2}]}\right\}.\end{aligned}$$

Set $q=1-(L+1)^{-1/2}$, $\gamma=1-\min(1/2,t^{1/4})$, and this, together with (8.6), verifies (8.5). $\square$

**Acknowledgments.** Suggestions of the Editor, Associate Editor and three referees are greatly appreciated.

Department of Mathematical Sciences
University of Texas at Dallas
Richardson, Texas 75083-0688
USA
E-mail: efrom@utdallas.edu